\documentclass{article}
\usepackage{amsmath,amsthm,graphicx,url}
 \twocolumn
\makeatletter
\long\def\proofbox#1{\gdef\@proofbox{#1}}

 \def\affil#1{\\{\small\sl#1\par}}
 \long\def\author#1{\gdef\@author{#1}}
 \author{Tommaso Toffoli ({\tt tt\char"40bu.edu})\affil{Electrical and
Computer Engineering, Boston University, Boston, MA 02215}}

 \long\def\abstract#1{\gdef\@abstract{#1}}
 \abstract{}

 \def\leftfrac{.32}
 \def\ritefrac{.68}

 \def\@maketitle{\newpage\noindent\leavevmode
  \begin{minipage}[t]{\leftfrac\textwidth}
    \hrule height0pt
    \@proofbox
  \end{minipage}\hfil
  \begin{minipage}[t]{\ritefrac\textwidth}
    \hrule height0pt
    \raggedleft
    \LARGE\@title\par
    \vskip4pt
    \large\@author
  \end{minipage}
  \vskip8pt
  \ifx\@abstract\@empty\else{\vskip.5em\leftskip1.25in\parskip4pt\small\@abstract\par\vskip.5em}\fi
  \noindent
  \rule{\textwidth}{0.4pt}
  \vskip16pt}
\makeatother

 \makeatletter

 \DeclareRobustCommand\em
        {\@nomath\em \ifdim \fontdimen\@ne\font >\z@
                       \upshape \else \slshape \fi}
 \newcommand{\ie}{i.e.,}

 \def\pages#1{}

 \newcommand{\sectlabel}[1]{\label{sect:#1}}
 \newcommand{\footlabel}[1]{\label{foot:#1}}
 \newcommand{\eqlabel}[1]{\label{eq:#1}}

 \newcommand{\Sect}[2][]{\def\t@mp{#1}%
\section{#2} \ifx\t@mp\@empty\else\sectlabel{#1}\fi}
 \newcommand{\Subsect}[2][]{\def\t@mp{#1}%
\subsection{#2} \ifx\t@mp\@empty\else\sectlabel{#1}\fi}
 \newcommand{\Foot}[2][]{\def\t@mp{#1}%
\footnote{#2} \ifx\t@mp\@empty\else\footlabel{#1}\fi}
 \newcommand{\Eq}[2][]{\def\t@mp{#1}%
\begin{equation}#2\ifx\t@mp\@empty\notag\else\eqlabel{#1}\fi\end{equation}}
 \newcommand{\Eqaligned}[2][]{\def\t@mp{#1}%
\begin{equation}\begin{aligned}#2\end{aligned}
\ifx\t@mp\@empty\notag\else\eqlabel{#1}\fi
\end{equation}}
 
\makeatother

 \def\pages#1{}
 \def\hunker{\lower1pt}

 \proofbox{{\small A somewhat abridged and adapted version of this review
appeared under the title ``Honesty in inference'' in {\sl American Scientist
\bf92}:2 (Feb.--Mar.\ 2004), 182--185.}}

 \title{Maxwell's daemon, the Turing machine,\\ and Jaynes' robot}
 \author{Tommaso Toffoli ({\tt tt\char"40bu.edu})\affil{Electr.\ and Computer
Eng., Boston University}}%

\begin{document}
\maketitle

\noindent {\bf Probability Theory---The Logic of Science}, E. T. Jaynes,
xxx+727 pp. Cambridge University Press, 2003. \$60.00.

\bigskip

 \def\myquote#1{\hfill\begin{minipage}{3.5in}\small#1\par\end{minipage}\bigskip}
 \myquote{{\sl The true Logic for this world is the Calculus of
Probabilities.} [Maxwell, to his friend Lewis Cambell]}

\noindent The review of a posthumous book (Edwin Jaynes died in 1998) demands
particular care in balancing love for the author with critical detachment in
the interest of the readers. I risk, of course, to be found wanting on both
counts.

Edwin Jaynes is best known as pioneer and champion of the Maximum Entropy
Principle ({\sc maxent} for insiders): If, in a given context, you need to
formulate a probability distribution on which to base your initial bets,
choose, among all possible distributions that agree with what you know about
the problem, the one having maximum entropy. Why? Is this guaranteed to be
the ``real'' (whatever that may mean) probability distribution?  Of course
not! In fact you will most likely replace it with a new one as soon as you
see the outcome of the first trial---because by then you'll have one more
piece of information. Why, then? Because any other choice would be
equivalent to \emph{throwing away} some of the information you have or
assuming information you \emph{don't} have. If you had to program a robot (of
which more below) to autonomously cope with the uncertainties of Mars's
environment, you'd want its choices to be weigthed by the {\sc maxent}
principle; anything else would be indefensible.

In view of the central role that this principle plays in Jaynes' gospel, it
may come as a surprise that, out of the book's 700 pages, barely 30 are
devoted to it; important as {\sc maxent} is in practice, from a conceptual
viewpoint Jaynes clearly regards it a mere corollary. What, then, is this book
about?  In Jaynes' words, ``simply, \emph{probability theory as extended
logic}.''

\bigskip

In introductory philosophy courses it is not unusual to see induction naively
presented as ``the opposite of'' deduction. Deduction always draws assured
consequences from the premises, but---under a curse akin to
Cassandra's---will never tell us anything that was not already in the
premises; specifically, if the premises represent what we know, deduction
cannot tell us anything we didn't know to begin with.  Induction, on the
other hand, works backwards from empirical facts to general
principles.  In this way it can tell us something new, but---that curse
again---will never guarantee its revelations.

As a matter of fact, \emph{pace} Hume and Popper, no such antithesis
exists. There is a single inference machine, whose kernel is \emph{ordinary
logic}. We feed the machine with a collection of statements---logic does not
care in what order they are given and whether we call them hypotheses,
empirical facts, axioms, or priors---and we ask a Yes/No question $A$. Then
we turn the crank, and after some grinding the machine will reply ``Yes,''
``No,'' or ``I don't know.'' (If you worry about undecidable problems, add a
countdown timer, and when the time is up the machine will answer ``I don't
know \emph{yet}.'')  \emph{That's all there is to it!}

Owing to the above curse, a single shot at this machine is of limited
value. Its glory is revealed, however, when we make it part of a feedback
loop.  Suppose we got ``Don't know'' for an answer. We can then further and
further qualify the input statements---the premises---until we get a definite
Yes or No answer.  Proceeding in this way with differents sets of
premises---different scenarios, as it were---we can explore a whole space of
scenarios and see which ones yield ``Yes'' and which ``No.'' And where do we
get probability?  you may ask. Briefly (and here I must skip a lot of fine
print), once we have \emph{many} scenarios (think of all possible hands in a
card game or molecular configurations of a liter of gas) we may have to give
up trying to handle them one by one and resign ourselves instead to dealing
with them ``in bulk,'' bagging them according to some criteria, labeling each
bag and keeping track only of \emph{how many} scenarios we put in
it. (Absolute counts don't matter much here---only relative ones do. For
instance, if we add to our parameters a binary variable irrelevant to our
inquiry, each scenario will split into two and thus the size of every bag
will double, but the overall fraction of scenarios represented by each bag
will not.) The \emph{probability} of our original question $A$ being true is
the fraction of scenarios for which the machine answers ``Yes.''  This
probability is of course conditioned by the \emph{prior}, that is, the
premises $X$ that remained fixed while we were running through the scenarios
(one of these premises being the set of admissible scenarios itself), and is
accordingly denoted $p(A|X)$. If this probability distribution $p(A|X)$
doesn't seem to lead us to good bets, all we can do is revise our prior
$X$. ``As Harold Jeffreys explained long ago, induction is most valuable to a
scientist just when it turns out to be wrong; only then do we get new
fundamental knowledge.''

The probability thus defined obeys \emph{by construction} those
quantitative rules that, as Jaynes argues following Cox and Polya, are
the basic desiderata of plausible reasoning, namely the \emph{product rule},
 \Eq{
	p(AB|X) = p(A|X)p(B|AX) = p(B|X)p(A|BX)
 }
 and the \emph{sum rule}
 \Eq{
	p(A+B|X) = p(A|X)+p(B|X)-p(AB|X).
 }

 Of course, different priors may lead to different probability distributions.
However, the great majority of variables in the whole world may be presumed
to be roughly orthogonal to the specific phenomena under investigation. To
paraphrase Borel, the microscopic scenarios that make up a bag marked, say,
``coin turned up heads'' are likely to be totally replaced by different
microscopic scenarios when the prior changes from ``full moon'' to ``new
moon,'' yet the numerical contents of each bag will hardly be affected---and
thus the typical gambler in his den will never notice the difference.

This observation is conceivably at the root of the \emph{mind projection
fallacy} most often attacked by Jaynes, namely, that probabilities (such as
of a coin's turning up heads) represent \emph{intrinsic properties of
physics} rather than \emph{a description of one's knowledge (or lack of
knowledge)} of the situation. ``As we've noted several times, the idea that
probabilities are physical real things, based ultimately on observed
frequencies of random variables, underlies most recent expositions of
probability theory, which would seem to make it a branch of experimental
science\dots. We take time off for an interlude of physical considerations
that show the fundamental difficulty with the notion of a `random'
experiment.'' And there follows an amusing chapter in which Jaynes discusses
how to cheat at coin and die tossing, dicusses card hands, and shows that the
frequentists' cherished notion of probability as the `limit of an infinite
sequence of identical random experiments' is immediately dumped by the
frequentists themselves as soon as they are confronted with a questionable
shuffler.

Are probabilities, then, just in our mind? ``Our probabilities and the
entropies based on them are indeed `subjective' in the sense that they
represent human information; if they did not, they could not serve their
purpose. But they are completely `objective' in the sense that they are
determined \emph{by the information specified} [\ie\ the prior $X$, via the
inference machine], independent of anyone's personality, opinions, or
hopes. It is `objectivity' in this sense that we need if information is ever
to be a basis for new theoretical developments in science.''

\medskip

Much as Maxwell's daemon was introduced to bring attention to certain issues
of thermodynamics, and Turing's machine to discuss thought without being
conditioned by too much human baggage, thus to wean us from the
subjective/objective dilemma Jaynes encourages us to come up with inference
rules that could be implemented by a \emph{robot}.  ``In order to direct
attention to constructive things and away from controversial irrelevancies,
we shall invent an imaginary being. Its brain is to be designed \emph{by us},
so that it reasons according to certain definite rules. These rules will be
deduced from simple desiderata \dots\ that a rational person [would not wish
to violate].'' These desiderata boil down to

(I) \emph{Degrees of plausibility are represented by real numbers}.

(II) \emph{Qualitative correspondence with common sense}.

(III) \emph{Consistency}.

Specifically (and here Asimov's ``Three Laws of Robotics'' will come to
mind),

(IIIb) \emph{The robot always takes into account all the evidence it has
relevant to a question. It does not arbitrarily ignore some of the
information.\ \dots\ In other words, the robot is completely nonideological.}

(IIIc) \emph{The robot always represents equivalent states of knowledge by
equivalent plausibility assignments\dots.}

From these desiderata Jaynes naturally derives all of the Bayesian machinery
for statistical inference (incidentally, Bayes never wrote ``Bayes'
theorem''---it's really Laplace's). 

\medskip

Was all this necessary? I'll quote Jaynes from one of his latest
apologies for his life's enterprise (``Probabilities in quantum theory,''
1990; Jaynes retired in 1992, and after that he doesn't seem to have been
very active in probability theory, though he spent much time on a book on
{\sl The Physical Basis of Music}, available online.)

 \begin{quote} The original view of James Bernoulli and Laplace was that
probability theory is an extension of logic to the case where, because of
incomplete information, deductive reasoning by the Aristotelean syllogism is
not available. It was sometimes called ``the calculus of inductive
reasoning.'' All of Lapalace's great contributions to science were made with
the help of probability theory interpreted in this way.

But, starting in the mid-19th Century, Laplace's viewpoint came under violent
attack from Leslie Ellis, John Venn, George Boole, R. von Mises,
R. A. Fisher, M. G. Kendall, W. Feller, J. Neyman, and others\dots. This
school of thought was so aggressive that it has dominated the field almost
totally, so that virtually all probability textbooks in current use are
written from a viewpoint which rejects Laplace's interpretation and tries to
deny us the use of his methods. Almost the only exceptions are found in the
works of Harold Jeffreys and Arnold Zellner. We have written two short
histories of these matters, engaged in a polemical debate on them, and are
trying to finish a two-volume treatise on the subject, entitled {\sl Probability
Theory---The Logic of Science}.
 \end{quote}
 Here is the book at last, a big tome in two parts if not two volumes, edited
from a coherent though incomplete manuscript (formerly available online) by
Jaynes' one-time doctoral advisee and then colleague and collaborator Larry
Bretthorst.

\bigskip

Even after having slashed at half the world, in the end Jaynes finds himself
in very good company. ``Our system of probability could hardly be more
different from that of Kolmogorov, in style, philosophy, and purpose. What we
consider to be fully half of probability theory \dots---the principles for
assigning probabilities by logical analysis of incomplete information---is
not present at all in the Kolmogorov system. Yet, when all is said and done,
we find ourselves, to our own surprise, in agreement with Kolmogorov and in
disagreement with his critics, on nearly all technical issues.''

\medskip

Jaynes has been accused of exaggerating the differences between the
``orthodox'' probability school---Fisher and frequentist company---and his
{\sc maxent} school of Bayesian/Laplacian revivalists. But now that it no
longer takes an act of courage to profess the Bayesian faith (younger
students of statistics don't even realize that there ever was an issue) it's
hard to imagine how different the world of statistical inference might have
been, especially in image processing, data mining, and bioengineering,
without Jaynes' thirty years of preaching and rallying the flock.

\bigskip

Beside the polemical aspects, whose raison-d'\^etre is to some extent
extinguished today, and the conceptual and technical aspects, which remain
valid but may no longer be as revolutionary as they sound simply because the
revolution has been won, there is one facet of Jaynes' gospel that is
particularly visible in this book and remains as fresh as ever, namely, a
relentless drive to demystify and democratize science. Science is especially
successful when it turns difficult tasks into trivial ones (how else could we
be made free to tackle today's new difficult tasks?).  But technical people
who used to thrive on the very former difficulties are loath to lose their
competitive advantage. Abstruseness, lack of universal rules
(``ad-hoc-ness''), and intimation of magic betray an attempt (inconscious,
perhaps) to retain control.  ``Feller's perception was so keen that in
virtually every problem he was able to see a clever trick; and then gave only
the clever trick. So his readers got the impression that: (1)~Probability
theory has no systematic method; it is a collection of unrelated clever
tricks, each of which works on one problem but not on the next one;
(2)~Feller was possessed of superhuman cleverness; (3)~Only a person with
such cleverness can hope to find new useful results in probability theory.''

We can also sense the indignation that a poor Iowa orphan boy must have felt
towards those Pooh-Bahs (``Sir'' Ronald Fisher, Cambridge alumnus, F.R.S,
Royal Medal, and so forth) that, instead of using their high position to help
streamline scientific thought, kept adding clutter to it (and, to add insult
to injury, had the gall to knowingly speak of Nature's ``propensities''
without being physicists themselves).

\bigskip

In this book as well as in his other writings Jaynes displays an intensity of
missionary zeal worthy of a St.\ Paul. We can only speculate what vision hit
the latter on his road to Damascus, but in Jaynes' case we are in possession
of an autobiographical snapshot. (In all that follows, italics are mine.)

``Those who cling to a belief in the existence of `physical probabilities'
may react to [the present approach to probability] by pointing to quantum
theory, in which physical probabilities appear to express the most
fundamental laws of physics\dots. [But here too] probabilities [must] express
our ignorance.'' And he goes on to argue that ``this ignorance may be
unavoidable in practice, but in our present state of knowledge we do not know
whether it is unavoidable \emph{in principle}'' (remark that Jaynes is well
aware of Bell's theorem and has no qualms with it); ``the `central dogma'
simply asserts this, and draws the conclusion that belief in causes, and
searching for them, is \emph{philosophically na\"\i ve}.'' (Intimidation by
the ``orthodox,'' again.)  ``Because of recent spectacular advances in the
technology of experimentation, with increasingly detailed control over the
initial state of individual atoms \dots\ we thing that the stage is going to
be set, before many more years have passed, for the same thing [as the
Bayesian revolution] to happen in quantum theory.''  (In fact, today we are
beginning to dabble with nanotechnology and quantum computers.) ``A century
from now the true cause of microphenomena will be known to every schoolboy
and, to paraphrase Seneca, they will be incredulous that such clear truths
could have escaped us throughout the 20th (and into the 21st) century.''

And here is the snapshot, from the above-mentioned apology: ``This conviction
\emph{has affected the whole course of my career.} I had intended originally
to specialize in Quantum Electrodynamics, but this proved to be impossible
[because] whenever I look at any quantum-mechanical calculations, the basic
craziness of what we are doing \emph{rises in my gorge}[!] and I have to try
to find some different way of looking at the problem that makes physical
sense. Gradually, I came to see that the foundation of probability theory and
the role of human information have to be brought in, and I have spent many
years trying to understand them in the greatest generality.''

\bigskip

Then the consolidation of Bayesian inference, {\sc maxent}, this book, and
all the rest, have just been a long detour on a (yet unfulfilled) quest for
the Holy Grail?  What then? Thanks to the boy from Iowa, today every
schoolboy (and every scholar) can approach inference unhampered by absurd
probability myths.

Most thick technical books are purchased with hope but then hardly ever
touched. This book was meant to be read and is a pleasure to read. (Also,
since the subject index is rather skimpy, often the best way to locate an
item is to read a whole chapter through, whereby you'll find lots of extra
rewarding stuff.) The bibliography is rich and full of annotations. There are
many exercises, and occasionally the editor has creatively turned a gap in
the manuscript into an exercise for the reader.

You may disagree here and there, but Jaynes is an honest scientist and a good
teacher. On both counts, isn't this what {\sc maxent} is all about?

 \end{document}